\documentclass[12pt]{article}

\usepackage[english]{babel}
\usepackage{amsmath, amsfonts, amssymb, amsthm}
\usepackage{url}
\usepackage[misc,geometry]{ifsym} 
\usepackage{enumitem}
\usepackage{fullpage}
\usepackage{mathtools}

\usepackage{natbib}

\usepackage{algorithm}
\usepackage{algorithmic}
\usepackage{caption}

\newcommand{\R}{{\mathbb R}}

\newcommand{\myproof}{{\it Proof}  }
\newcommand{\e}{\varepsilon}
\newcommand{\la}{\langle}
\newcommand{\ra}{\rangle}
\def\Ec{\mathcal E}
\def\Hc{\mathcal H}

\def\tf{\tilde{f}}
\def\tg{\tilde{g}}

\newtheorem{Lm}{Lemma}[section]

\newtheorem{Th}{Theorem}[section]

\newtheorem{Cor}{Corollary}[section]

\theoremstyle{definition}
\newtheorem{Def}{Definition}[section]

\title{Gradient Method With Inexact Oracle for Composite Non-Convex Optimization} 

\author{Pavel Dvurechensky
				\thanks{Weierstrass Institute for Applied Analysis and Stochastics, Berlin; Institute for Information Transmission Problems RAS, 
								Moscow, pavel.dvurechensky@wias-berlin.de } 
				}

\date{\today}

\begin{document}

\maketitle

\begin{abstract}
In this paper, we develop new first-order method for composite non-convex minimization problems with simple constraints and inexact oracle. The objective function is given as a sum of "`hard"', possibly non-convex part, and "`simple"' convex part.
Informally speaking, oracle inexactness means that, for the "`hard"' part, at any point we can approximately calculate the value of the function and construct a quadratic function, which approximately bounds this function from above. 
We give several examples of such inexactness: smooth non-convex functions with inexact H\"older-continuous gradient, functions given by auxiliary uniformly concave maximization problem, which can be solved only approximately.
For the introduced class of problems, we propose a gradient-type method, which allows to use different proximal setup to adapt to geometry of the feasible set, adaptively chooses controlled oracle error, allows for inexact proximal mapping. 
We provide convergence rate for our method in terms of the norm of generalized gradient mapping and show that, in the case of inexact H\"older-continuous gradient, our method is universal with respect to H\"older parameters of the problem. Finally, in a particular case, we show that small value of the norm of generalized gradient mapping at a point means that a necessary condition of local minimum approximately holds at that point.
\end{abstract}


\section*{Introduction}
\label{S:Intro}
In this paper, we introduce new first-order method for non-convex composite optimization problems with inexact oracle. 
Namely, our problem of interest is as follows
\begin{equation}
\min_{x \in X \subseteq \Ec} \{ \psi(x) := f(x) + h(x)\},
\label{eq:PrStateInit}
\end{equation}
where $X$ is a closed convex set, $h(x)$ is a simple convex function, e.g. $\|x\|_1$. We assume that $f(x)$ is a general function endowed with an inexact first-order oracle, which is defined below (see Definition \ref{D:IO}). Informally speaking, at any point we can approximately calculate the value of the function and construct a quadratic function, which approximately bounds our $f(x)$ from above. An example of problem with this kind of inexactness is given in \cite{bogolubsky2016learning}, where the authors study a learning problem for parametric PageRank model.
 
First-order methods are widely developed since the earliest years of optimization theory, see, e.g., \cite{polyak1963gradient}. Recent renaissance in their development started more than ten years ago and was mostly motivated by fast growing problem sizes in applications such as Machine Learning, Data Analysis, Telecommunications. For many years, researchers mostly considered convex optimization problems since they have good structure and allow to estimate rate of convergence for proposed algorithms. Recently, non-convex problems started to attract fast growing attention, as they appear often in Machine Learning, especially in Deep Learning. Thus, high standards of research on algorithms for convex optimization started to influence non-convex optimization. Namely, it have become very important for newly developed methods to obtain a rate of convergence with respect to some criterion. Usually, this criterion is the norm of gradient mapping, which is a generalization of gradient for constrained problems, see, e.g. \cite{nesterov2004introduction}.

Already in \cite{polyak1987introduction}, the author analyzed how different types of inexactness in gradient values influence gradient method for unconstrained smooth convex problems.
At the moment, theory for convex optimization algorithms with inexact oracle is well-developed in a series of papers \cite{aspremont2008smooth,devolder2014first,dvurechensky2016stochastic}. In \cite{aspremont2008smooth}, it was proposed to calculate inexactly the gradient of the objective function and extend Fast Gradient Method of \cite{nesterov2005smooth} to be able to use inexact oracle information. In \cite{devolder2014first}, a general concept of inexact oracle is introduced for convex problems, Primal, Dual and Fast gradient methods are analyzed. In \cite{dvurechensky2016stochastic}, the authors develop Stochastic Intermediate Gradient Method for problems with stochastic inexact oracle, which provides good flexibility for solving convex and strongly convex problems with both deterministic and stochastic inexactness. 

The theory for non-convex smooth, non-smooth and stochastic problems is well developed in \cite{ghadimi2016accelerated,ghadimi2016mini-batch}. In \cite{ghadimi2016accelerated}, problems of the form \eqref{eq:PrStateInit}, where $X \equiv \R^n$ and $f(x)$ is a smooth non-convex function are considered in the case when the gradient of $f(x)$ is exactly available, as well as when it is available through stochastic approximation. Later, in \cite{ghadimi2016mini-batch} the authors generalized these methods for constrained problems of the form \eqref{eq:PrStateInit} in both deterministic and stochastic settings. 

Nevertheless, it seems to us that gradient methods for non-convex optimization problems with deterministic inexact oracle lack sufficient development. The goal of this paper is to fill this gap.

It turns out that smooth minimization with inexact oracle is closely connected with minimization of functions with H\"older-continuous gradient. 
We say that a function $f(x)$ has H\"older-continuous gradient on $X$ iff there exist $\nu \in [0,1]$ and $L_{\nu} \geq 0$ s.t.
$$
\|\nabla f(x) - \nabla f(y)\|_{\Ec,*} \leq L_{\nu}\|x-y\|_{\Ec}^{\nu}, \quad x,y \in X.
$$
In \cite{devolder2014first} it was shown that a convex problem with H\"older-continuous subgradient can be considered as a smooth problem with deterministic inexact oracle. Later, universal gradient methods for convex problems with H\"older-continuous subgradient were proposed in \cite{nesterov2015universal}. These algorithms do not require to know H\"older parameter $\nu$ and H\"older constant $L_{\nu}$. Thus, they are universal with respect to these parameters. \cite{lan2015generalized} proposed methods for non-convex problems of the form \eqref{eq:PrStateInit}, where $f(x)$ has H\"older-continuous gradient. These methods rely on Euclidean norm and are good when the euclidean projection onto the set $X$ is simple.

Our contribution in this paper is as follows.
\begin{enumerate}
	\item We generalize for non-convex case the definition of inexact oracle in \cite{devolder2014first} and provide several examples, where such inexactness can arise. We consider two types of errors -- controlled errors, which can be made as small as desired, and uncontrolled errors, which can only be estimated.
	\item We introduce new gradient method for problem \eqref{eq:PrStateInit} and prove a theorem (see Theorem \ref{Th:UGMRate}) on its rate of convergence in terms of the norm of generalized gradient mapping. Our method is adaptive to the controlled oracle error, is capable to work with inexact proximal mapping, has flexibility of choice of proximal setup, based on the geometry of set $X$. 
	\item We show that, in the case of problems with inexact H\"older-continuous gradient, our method is universal, that is, it does not require to know in advance a H\"older parameter $\nu$ and H\"older constant $L_{\nu}$ for the function $f(x)$, but provides best known convergence rate uniformly in H\"older parameter $\nu$.  
\end{enumerate}
Thus, we provide a universal algorithm for non-convex H\"older-smooth composite optimization problems with deterministic inexact oracle.

The rest of the paper is organized as follows. In Section \ref{S:IO}, we define deterministic inexact oracle for non-convex problems and provide several examples. In Section \ref{S:Alg}, we describe our algorithm, prove the convergence theorem. Also we provide two corollaries for particular cases of smooth functions and H\"older-smooth functions. Note that the latter case includes the former one. Finally, we provide some explanations about how convergence of the norm of generalized gradient mapping to zero leads to a good approximation for a point, where a necessary optimality condition for Problem \eqref{eq:PrStateInit} holds. Note that we use different reasoning from what can be found in literature. 

\textbf{Notation}
Let $\Ec$ be a finite-dimensional real vector space and $\Ec^*$ be its dual. We denote the value of linear function $g \in \Ec^*$ at $x\in \Ec$ by $\la g, x \ra$. Let $\|\cdot\|_{\Ec}$ be some norm on $\Ec$, $\|\cdot\|_{\Ec,*}$ be its dual.

\section{Inexact Oracle}
\label{S:IO}
In this section, we define the inexact oracle and describe several examples where it naturally arises.
\begin{Def}
\label{D:IO}
We say that a function $f(x)$ is equipped with an \textit{inexact first-order oracle} on a set $X$ if there exists $\delta_u > 0$ and at any point $x \in X$ for any number $\delta_c > 0$ there exists a constant $L(\delta_c) \in (0, +\infty)$ and one can calculate  $\tf(x,\delta_c,\delta_u) \in \R$ and $\tg(x,\delta_c,\delta_u) \in \Ec^*$ satisfying
\begin{align}
	&|f(x) - \tf(x,\delta_c,\delta_u)| \leq \delta_c+\delta_u, 
	\label{eq:dL_or_def_2} \\
	&f(y)-(\tf(x,\delta_c,\delta_u)  - \la\tg(x,\delta_c,\delta_u) ,y-x \ra) \leq \frac{L(\delta_c)}{2}\|x-y\|_{\Ec}^2 + \delta_c+\delta_u, \quad \forall y \in X.
  \label{eq:dL_or_def_1}
\end{align}
\end{Def}
In this definition, $\delta_c$ represents the error of the oracle, which we can control and make as small as we would like to. On the opposite, $\delta_u$ represents the error, which we can not control. 
The idea behind the definition is that at any point we can approximately calculate the value of the function and construct an upper quadratic bound.


Let us now consider several examples.

\subsection{Smooth Function with Inexact Oracle Values}
\label{SS:SFIO}
Let us assume that 
\begin{enumerate}
	\item Function $f(x)$ is $L$-smooth on $X$, i.e. it is differentiable and, for all $x,y \in X$, $\|\nabla f(x) - \nabla f(y)\|_{{\Ec},*} \leq L\|x-y\|_{\Ec}$.
	\item Set $X$ is bounded with $\max_{x,y \in X} \|x-y\|_{\Ec} \leq D$.
	\item There exist $\bar{\delta}_u^1, \bar{\delta}_u^2 > 0$ and at any point $x \in Q$, for any $\bar{\delta}_c^1, \bar{\delta}_c^2 > 0$, we can calculate approximations $\bar{f}(x)$ and $\bar{g}(x)$ s.t. $|\bar{f}(x) - f(x)| \leq \bar{\delta}_c^1+\bar{\delta}_u^1$, $\|\bar{g}(x)-\nabla f(x)\|_{{\Ec},*} \leq \bar{\delta}_c^2+\bar{\delta}_u^2$.
\end{enumerate} 

Then, using  $L$-smoothness of $f(x)$, we obtain, for any $y \in X$, 
\begin{align}
f(y) &\leq f(x) + \la \nabla f(x), y-x \ra + \frac{L}{2}\|x-y\|_{\Ec}^2 \\
     &\leq \bar{f}(x) + \bar{\delta}_c^1+\bar{\delta}_u^1 + \la \nabla \bar{g}(x), y-x \ra + \la \nabla f(x) - \bar{g}(x), y-x \ra + \frac{L}{2}\|x-y\|_{\Ec}^2 \\
	   &\leq \bar{f}(x)  + \la \nabla \bar{g}(x), y-x \ra + \frac{L}{2}\|x-y\|_{\Ec}^2 + \bar{\delta}_c^1+\bar{\delta}_u^1 + (\bar{\delta}_c^2+\bar{\delta}_u^2 ) D.
\end{align}
Thus, $(\bar{f}(x),\bar{g}(x))$ is an inexact first-order oracle with $\delta_u = \bar{\delta}_u^1+\bar{\delta}_u^2D$,  $\delta_c = \bar{\delta}_c^1+\bar{\delta}_c^2D$, and $L(\delta_c) \equiv L$.

\subsection{Smooth Function with H\"older-Continuous Gradient}
\label{SS:HCG}
Assume that $f(x)$ is differentiable and its gradient is H\"older-continuous, i.e. for some $\nu \in [0,1]$  and $L_{\nu} \geq 0$,
\begin{equation}
\|\nabla f(x) - \nabla f(y) \|_* \leq L_{\nu} \|x-y\|_{\Ec}^{\nu}, \forall x,y \in X.
\label{eq:nfH}
\end{equation} 
Then 
\begin{equation}
f(y) \leq f(x)  + \la \nabla f(x) ,y-x \ra + \frac{L_{\nu}}{1+\nu}\|x-y\|_{\Ec}^{1+\nu}, \quad \forall x,y \in X.
\label{eq:nfH1}
\end{equation}
It can be shown, see \cite{nesterov2015universal}, Lemma 2, that, for all $x \in X$ and any $\delta >0$,
\begin{equation}
f(y)-(f(x)  - \la \nabla f(x) ,y-x \ra) \leq \frac{L(\delta)}{2}\|x-y\|_{\Ec}^2 + \delta, \quad \forall y \in X,
\label{eq:nfH2}
\end{equation}
where 
\begin{equation}
L(\delta) = \left( \frac{1-\nu}{1+\nu} \cdot \frac{2}{\delta} \right)^{\frac{1-\nu}{1+\nu}} L_{\nu}^{\frac{2}{1+\nu}}.
\label{eq:Lofd}
\end{equation}
Thus, $(f(x),\nabla f(x))$ is an inexact first-order oracle with $\delta_u = 0$,  $\delta_c = \delta$, and $L(\delta)$ given by \eqref{eq:Lofd}.

Note that, if $(f(x),\nabla f(x))$ can only be calculated inexactly as in Subsection \ref{SS:SFIO}, their approximations will again be an inexact first-order oracle.
%

\subsection{Function Given by Maximization Subproblem}
\label{SS:Sm}
Assume that function $f(x): \Ec \to \R$ is defined by an auxiliary optimization problem
\begin{equation}
f(x) = \max_{u \in U \subseteq \Hc} \{\Psi(x,u):= -G(u) + \la A u , x \ra \}, 
\label{eq:fsmtfr}
\end{equation}
where $A: \Hc \to \Ec^*$ is a linear operator, $G: \Hc \to R$ is a continuously differentiable uniformly convex function of degree $\rho \geq 2$ with parameter $\sigma_{\rho} \geq 0$. The last means that 
\begin{equation}
\la \nabla G(u_1) - \nabla G(u_2) , u_1 - u_2 \ra \geq \sigma_{\rho} \|u_1-u_2\|_{\Hc}^{\rho}, \quad \forall u_1, u_2 \in U,
\label{eq:ucdef}
\end{equation}
where $\|\cdot\|_{\Hc}$ is some norm on $\Hc$. Note that $f(x)$ is differentiable and $\nabla f(x) = A u^*(x)$, where $u^*(x)$ is the optimal solution in \eqref{eq:fsmtfr} for fixed $x$. 

Extending the proof in \cite{nesterov2015universal}, we can prove the following.
\begin{Lm}
If $G$ is uniformly convex on $X$, then the gradient of $f$ is H\"older-continuous with
\begin{equation}
\nu = \frac{1}{\rho-1}, \quad L_{\nu} = \frac{\|A\|_{\Hc \to \Ec^*}^{\frac{\rho}{\rho-1}}}{\sigma_{\rho}^{\frac{1}{\rho-1}}},
\label{eq:nfHc}
\end{equation}
where $\|A\|_{\Hc \to \Ec^*} = \max \{\|Au\|_{\Ec,*}: \|u\|_{\Hc}=1\}$.
\label{Lm:nfH}
\end{Lm}
\myproof
From the optimality conditions in \eqref{eq:fsmtfr}, we obtain
\begin{align}
&\la A^T x_1 - \nabla G(u(x_1)), u(x_2) - u(x_1) \ra \leq 0, \\
&\la A^T x_2 - \nabla G(u(x_2)), u(x_1) - u(x_2) \ra \leq 0.
\end{align}
Adding these inequalities, we obtain, by definition of uniformly convex function,
\begin{equation}
\la A^T (x_1-x_2), u(x_1) - u(x_2)  \ra \geq \la \nabla G(u(x_1)) - \nabla G(u(x_2)), u(x_1) - u(x_2)  \ra \stackrel{\eqref{eq:ucdef}}{\geq} \sigma_{\rho} \|u(x_1) - u(x_2)\|_{\Hc}^{\rho}.
\end{equation}
on the other hand,
\begin{align}
\|A(u(x_1) - u(x_2))\|_{\Ec,*}^2 
&\leq \|A\|_{\Hc \to \Ec^*}^2 \|u(x_1) - u(x_2)\|_{\Hc}^2 \\
&\leq \|A\|_{\Hc \to \Ec^*}^2 \left(\frac{1}{\sigma_{\rho}} \la A^T (x_1-x_2), u(x_1) - u(x_2)  \ra \right)^{2/\rho} \\
&\leq \frac{\|A\|_{\Hc \to \Ec^*}^2}{\sigma_{\rho}^{2/\rho}} \|A(u(x_1) - u(x_2))\|_{\Ec,*}^{2/\rho} \|x_1-x_2\|_{\Ec}^{2/\rho}.
\end{align}
Thus,
\begin{equation}
\|A(u(x_1) - u(x_2))\|_{\Ec,*}^{2-2/\rho} \leq  \frac{\|A\|_{\Hc \to \Ec^*}^2}{\sigma_{\rho}^{2/\rho}} \|x_1-x_2\|_{\Ec}^{2/\rho},
\end{equation}
which proves the Lemma.
\qed

Let us now consider a situation, when the maximization problem in \eqref{eq:fsmtfr} can be solved only inexactly by some auxiliary numerical method. It is natural to assume that, for any $x \in X$ and any $\delta > 0$, we can calculate a point $u_x \in U$ s.t. 
\begin{equation}
0 \leq f(x)-\Psi(x,u_x) = \Psi(x,u^*(x)) -\Psi(x,u_x) \leq \delta.
\label{eq:Psier}
\end{equation}
Since $\ln(t) $ is a concave function, for any $\rho \geq 2$ and $t,\tau \geq 0$, we have
\begin{equation}
\ln \left( \frac{1}{\rho}t^\rho+ \frac{\rho-1}{\rho} \tau^{\frac{\rho}{\rho-1}} \right) \geq \frac{1}{\rho} \ln \left( t^\rho\right) + \frac{\rho-1}{\rho} \ln\left(\tau^{\frac{\rho}{\rho-1}}\right) = \ln(t \tau).
\end{equation}
Using this inequality with
\begin{equation}
t=\sigma_{\rho}^{1/\rho} \|u^*(x)-u_x\|_{\Hc}, \quad \tau = \frac{\|A\|_{\Hc \to \Ec^*}}{\sigma_{\rho}^{1/\rho}} \|y-x\|_{\Ec},
\end{equation}
we obtain, for any $y \in X$,
\begin{align}
\la A(u^*(x)-u_x), y-x \ra &\leq \|A\|_{\Hc \to \Ec^*}\|u^*(x)-u_x\|_{\Hc} \|y-x\|_{\Ec} \\
& \leq \frac{\sigma_{\rho}}{\rho}\|u^*(x)-u_x\|_{\Hc}^{\rho} + \frac{\|A\|_{\Hc \to \Ec^*}^{\frac{\rho}{\rho-1}}}{\frac{\rho}{\rho-1}\sigma_{\rho}^{\frac{1}{\rho-1}}} \|y-x\|_{\Ec}^{\frac{\rho}{\rho-1}} \\
& = \frac{\sigma_{\rho}}{\rho}\|u^*(x)-u_x\|_{\Hc}^{\rho} + \frac{L_{\nu}}{1+\nu}\|y-x\|_{\Ec}^{1+\nu},
\end{align}
where $\nu$ and $L_{\nu}$ are defined in \eqref{eq:nfHc}. At the same time, since $\Psi(x,u)$ \eqref{eq:fsmtfr} is uniformly concave in second argument, we have
\begin{equation}
\frac{\sigma_{\rho}}{\rho}\|u^*(x)-u_x\|_{\Hc}^{\rho} \leq \Psi(x,u^*(x)) - \Psi(x, u_x) \stackrel{\eqref{eq:Psier}}{\leq} \delta.
\end{equation}
Combining this inequality with the previous one, we obtain
\begin{equation}
\la A(u^*(x)-u_x), y-x \ra \leq \frac{L_{\nu}}{1+\nu}\|x-y\|_{\Ec}^{1+\nu} + \delta.
\label{eq:Auxbound}
\end{equation}
Since $f$ has H\"older-continuous gradient with parameters \eqref{eq:nfHc}, using \eqref{eq:nfH1}, we obtain
\begin{align}
f(y) &\leq f(x) + \la \nabla f(x) ,y-x \ra + \frac{L_{\nu}}{1+\nu}\|x-y\|_{\Ec}^{1+\nu} \\
& \stackrel{\eqref{eq:Psier}}{\leq} \Psi(x,u_x) + \delta +  \la Au_x ,y-x \ra + \la A(u^*(x)-u_x), y-x \ra + \frac{2L_{\nu}}{1+\nu}\|x-y\|_{\Ec}^{1+\nu} \\
& \stackrel{\eqref{eq:Auxbound}}{\leq} \Psi(x,u_x) + \la Au_x ,y-x \ra + \frac{2L_{\nu}}{1+\nu}\|x-y\|_{\Ec}^{1+\nu}  + 2 \delta \\
& \stackrel{\eqref{eq:nfH1},\eqref{eq:nfH2},\eqref{eq:Lofd}}{\leq} \Psi(x,u_x) + \la Au_x ,y-x \ra + \frac{2L(\delta)}{2}\|x-y\|_{\Ec}^{2}  + 4 \delta.
\end{align}
Thus, we have obtained that $(\Psi(x,u_x), Au_x)$ is an inexact first-order oracle with $\delta_u = 0$,  $\delta_c = 4\delta$, and $L(\delta_c)$ given by \eqref{eq:Lofd} with $\delta = \delta_c/4$.

\section{Adaptive Gradient Method for Problems with Inexact Oracle}
\label{S:Alg}
To construct our algorithm for problem \eqref{eq:PrStateInit}, we introduce, as it usually done, proximal setup \cite{ben-tal2015lectures}.
We choose a {\it prox-function} $d(x)$ which is continuous, convex on $X$ and
\begin{enumerate}
	\item admits a continuous in $x \in X^0$ selection of subgradients 	$d'(x)$, where $x \in X^0 \subseteq X$  is the set of all $x$, where $d'(x)$ exists;
	\item $d(x)$ is $1$-strongly convex on $X$ with respect to $\|\cdot\|_{\Ec}$, i.e., for any $x \in X^0, y \in X$ $d(y)-d(x) -\la d'(x) ,y-x \ra \geq \frac12\|y-x\|_{\Ec}^2$.
\end{enumerate} 
We define also the corresponding {\it Bregman divergence} $V[z] (x) = d(x) - d(z) - \la d'(z), x - z \ra$, $x \in X, z \in X^0$. Standard proximal setups, i.e. Euclidean, entropy, $\ell_1/\ell_2$, simplex ,  nuclear norm, spectahedron can be found in \cite{ben-tal2015lectures}.
We will use Bregman divergence in so called {\it composite prox-mapping}
\begin{equation}
\min_{x \in X} \left\{\la g,x \ra + \frac{1}{\gamma} V[\bar{x}](x) +h(x) 		 \right\} ,
\label{eq:PrMap}
\end{equation}
where $\gamma >0 $, $\bar{x} \in X^0$, $g \in \Ec^*$ are given. We allow this problem to be solved inexactly in the following sense.
\begin{Def}
\label{D:IPM}
Assume that we are given $\delta_{pu} >0$, $\gamma >0 $, $\bar{x} \in X^0$, $g \in \Ec^*$. We call a point $\tilde{x} = \tilde{x}(\bar{x},g,\gamma,\delta_{pc},\delta_{pu}) \in X^0$ an {\it inexact composite prox-mapping} iff for any $\delta_{pc} >0$ we can calculate $\tilde{x}$ and there exists $p \in \partial h(\tilde{x})$ s.t. it holds that
\begin{equation}
\left\la g + \frac{1}{\gamma}\left[d'(\tilde{x}) - d'(\bar{x}) \right] + p, u - \tilde{x} \right\ra \geq - \delta_{pc}-\delta_{pu}, \quad \forall u \in X.
\label{eq:InPrMap}
\end{equation}
We write
\begin{equation}
\tilde{x} =  {\mathop {\arg \min }\limits_{x\in X}}^{\delta_{pc}+\delta_{pu}}\left\{\la g,x \ra + \frac{1}{\gamma} V[\bar{x}](x) +h(x) 		 \right\}
\label{eq:InPrMap1}
\end{equation}
and define
\begin{equation}
g_X (\bar{x},g,\gamma,\delta_{pc},\delta_{pu}) := \frac{1}{\gamma}(\bar{x}-\tilde{x}).
\label{eq:g_Q}
\end{equation}
\end{Def}
This is a generalization of inexact composite prox-mapping in \cite{ben-tal2015lectures}. Note that if $\tilde{x}$ is an exact solution of \eqref{eq:PrMap}, inequality \eqref{eq:InPrMap} holds with $\delta_{pc}=\delta_{pu}=0$ due to first-order optimality condition. Similarly to Definition \ref{D:IO}, $\delta_{pc}$ represents an error, which can be controlled and made as small as it is desired, $\delta_{pu}$ represents an error which can not be controlled.

Our main scheme is Algorithm \ref{Alg:UGM}. 

\begin{algorithm}[h!]
\caption{Adaptive Gradient Method for Problems with Inexact Oracle}
\label{Alg:UGM}
\begin{algorithmic}[1]
   \REQUIRE accuracy $\e > 0$, uncontrolled oracle error $\delta_u >0$, uncontrolled error of composite prox-mapping $\delta_{pu} >0$, starting point $x_0 \in X^0$, initial guess $L_0 >0$, prox-setup: $d(x)$ -- $1$-strongly convex w.r.t. $\|\cdot\|_{\Ec}$, $V[z] (x) := d(x) - d(z) - \la d'(z), x - z \ra$.
   \STATE Set $k=0$.
   \REPEAT
			\STATE Set $M_k=L_k/2$.
			\REPEAT
				\STATE Set $M_k=2M_k$, $\delta_{c,k} = \delta_{pc,k} = \frac{\e}{20M_k}$.
				\STATE Calculate $\tf(x_k,\delta_{c,k},\delta_u)$ and $\tg(x_k,\delta_{c,k},\delta_u)$.	
				\STATE Calculate 
				\begin{equation}
				w_k={\mathop {\arg \min }\limits_{x\in X}}^{\delta_{pc,k}+\delta_{pu}} \left\{\la \tg(x_k,\delta_{c,k},\delta_u),x \ra + M_kV[x_k](x) +h(x) 		 \right\}.
				\label{eq:UGMwStep}
				\end{equation}
				\STATE Calculate $\tf(w_k,\delta_{c,k},\delta_u)$. 
			\UNTIL{
			\begin{equation}
			\tf(w_k,\delta_{c,k},\delta_u) \leq \tf(x_k,\delta_{c,k},\delta_u) + \la \tg(x_k,\delta_{c,k},\delta_u) ,w_k - x_k \ra  +\frac{M_k}{2}\|w_k - x_k\|_{\Ec}^2 +	\frac{\e}{10M_k} + 2 \delta_u.
			\label{eq:UGMCheck}
			\end{equation}}
			\STATE Set $x_{k+1} = w_k$, $L_{k+1}=M_k/2$, $k=k+1$.
   \UNTIL{$\min_{i\in 1,...,k} \left\|M_i(x_i-x_{i+1})\right\|_{\Ec}   \leq \e $}
		\ENSURE The point $x_{K+1}$ s.t. $K = \arg \min_{i\in 1,...,k} \left\|M_i(x_i-x_{i+1})\right\|_{\Ec}$.	
\end{algorithmic}
\end{algorithm}

We will need the following simple extension of Lemma 1 in \cite{ghadimi2016mini-batch} to perform the theoretical analysis of our algorithm.
\begin{Lm}
\label{Lm:gr_map_pr_1}
Let $\tilde{x} = \tilde{x}(\bar{x},g,\gamma,\delta_{pc},\delta_{pu})$ be an inexact composite prox-mapping and $g_X (\bar{x},g,\gamma,\delta_{pc},\delta_{pu})$ be defined in \eqref{eq:g_Q}. Then, for any $\bar{x} \in X^0$, $g \in \Ec^*$ and $\gamma, \delta_{pc},\delta_{pu} > 0$, it holds
\begin{equation}
\gamma \la g, g_X (\bar{x},g,\gamma,\delta_{pc},\delta_{pu}) \ra \geq \gamma \|g_X (\bar{x},g,\gamma,\delta_{pc},\delta_{pu})\|_{\Ec}^2 + (h(\tilde{x}(\bar{x},g,\gamma,\delta_{pc},\delta_{pu}))-h(x)) - \delta_{pc}-\delta_{pu} .
\label{eq:gr_map_pr_1}
\end{equation}
\end{Lm}
\myproof
Taking $u=\bar{x}$ in \eqref{eq:InPrMap} and rearranging terms, we obtain, by convexity of $h(x)$ and strong convexity of $d(x)$,
\begin{align}
\la g, \bar{x}- \tilde{x} \ra & \geq \frac{1}{\gamma} \la d'(\tilde{x}) - d'(\bar{x}) , \tilde{x}  - \bar{x} \ra + \la p, \tilde{x}  - \bar{x} \ra - \delta_{pc}-\delta_{pu} \notag \\
& \geq \frac{1}{\gamma} \|\tilde{x}  - \bar{x}\|_{\Ec}^2 + (h(\tilde{x}) - h(\bar{x}))- \delta_{pc}-\delta_{pu}.
\end{align}
Applying the definition \eqref{eq:g_Q}, we finish the proof.
\qed

Now we state the main 
\begin{Th}
\label{Th:UGMRate}
Assume that $f(x)$ is equipped with an inexact first-order oracle in the sense of Definition \ref{D:IO} and for any constants $c_1, c_2 >0$ there exists an integer $i \geq 0$ s.t. $2^ic_1 \geq L\left(\frac{c_2}{c_12^i}\right)$. Assume also that there exists a number $\psi^* > -\infty$ such that $\psi(x) \geq \psi^*$ for all $x \in X$. Then, after $N$ iterations of Algorithm \ref{Alg:UGM}, it holds that
\begin{equation}
\left\|M_K (x_K - x_{K+1}))\right\|_{\Ec}^2  \leq \left(\sum_{k=0}^{N-1} \frac{1}{2M_k} \right)^{-1} (\psi(x_0) - \psi^* + N(4\delta_u + \delta_{pu})) + \frac{\e}{2}.
\label{eq:pg_la_2_rate}
\end{equation}
Moreover, the total number of checks of Inequality \eqref{eq:UGMCheck} is not more than 
\begin{equation}
2N-1+\log_2\frac{M_{N-1}}{L_0}.
\label{eq:InnerChecksBound}
\end{equation}
\end{Th}
\myproof
First of all let us show that the procedure of search of point $w_k$ satisfying \eqref{eq:UGMwStep}, \eqref{eq:UGMCheck} is finite. Let $i_k \geq 0$ be the current number of performed checks of inequality \eqref{eq:UGMCheck} on the step $k$. Then $M_k = 2^{i_k}L_k$. At the same time, by Definition \ref{D:IO} $L(\delta_{c,k}) = L\left( \frac{\e}{16M_k} \right) = L\left( \frac{\e}{16 \cdot 2^{i_k}L_k} \right)$. Hence, by the Theorem assumptions, there exists $i_k \geq 0$ s.t. $M_k = 2^{i_k}L_k \geq L(\delta_{c,k})$. At the same time, we have 
\begin{align}
\tf(w_k,\delta_{c,k},\delta_u) - \frac{\e}{20M_k} - \delta_u & \stackrel{\eqref{eq:dL_or_def_2}}{\leq} f(w_k) \\
&\stackrel{\eqref{eq:dL_or_def_1}}{\leq} \tf(x_k,\delta_{c,k},\delta_u) + \la \tg(x_k,\delta_{c,k},\delta_u),w_k - x_k \ra \\
& + \frac{L(\delta_{c,k})}{2}\|w_k - x_k\|_{\Ec}^2 + \frac{\e}{20M_k} + \delta_u,
\notag
\end{align}
which leads to \eqref{eq:UGMCheck} when $M_k \geq L(\delta_{c,k})$.

Let us now obtain the rate of convergence. We denote, for simplicity, $\tf_k = \tf(x_k,\delta_{c,k},\delta_u)$, $\tg_k = \tg(x_k,\delta_{c,k},\delta_u)$, $\tg_{X,k} =  g_X \left(x_k, \tg_k ,\frac{1}{M_k},\delta_{pc,k},\delta_{pu}\right)$
Note that
\begin{equation}
\tg_{X,k} \stackrel{\eqref{eq:InPrMap1},\eqref{eq:g_Q},\eqref{eq:UGMwStep}} {=} M_k (x_k - x_{k+1}).
\label{eq:gXeqMxmx}
\end{equation}  Using definition of $x_{k+1}$, we obtain, for any $k=0,\dots,N-1$,
\begin{align}
f(x_{k+1}) - \frac{\e}{20M_k} -\delta_u & = f(w_k) - \frac{\e}{20M_k} -\delta_u \\
&  \stackrel{\eqref{eq:dL_or_def_2}}{\leq} \tf(w_k,\delta_{c,k},\delta_u) \\
& \stackrel{\eqref{eq:UGMCheck}}{\leq} \tf_k +  \la \tg_k , x_{k+1}-x_k \ra +\frac{M_k}{2} \|x_{k+1}-x_k\|_{\Ec}^2 + \frac{\e}{10M_k} + 2\delta_u \\
& \stackrel{\eqref{eq:gXeqMxmx}}{=}   \tf_k   - \frac{1}{M_k} \left\la \tg_k  , \tg_{X,k} \right\ra +\frac{1}{2M_k} \left\|\tg_{X,k}\right\|_{\Ec}^2 + \frac{\e}{10M_k} + 2\delta_u \\
& \stackrel{\eqref{eq:dL_or_def_2},\eqref{eq:gr_map_pr_1}}{\leq}   f(x_k) + \frac{\e}{20M_k} + \delta_u - \left[\frac{1}{M_k} \left\|\tg_{X,k}\right\|_{\Ec}^2  + h(x_{k+1})-h(x_k) - \frac{\e}{20M_k} - \delta_{pu}\right]  \\
& + \frac{1}{2M_k} \left\|\tg_{X,k}\right\|_{\Ec}^2 + \frac{\e}{10M_k} + 2\delta_u \notag.
\end{align}
This leads to
$$
\psi(x_{k+1}) \leq  \psi(x_k) - \frac{1}{2M_k} \left\|\tg_{X,k}\right\|_{\Ec}^2 + \frac{\e}{4M_k} + 4\delta_u + \delta_{pu}, \quad k=0,\dots,N-1.
$$
Summing up these inequalities, we get
\begin{align}
&\left\|\tg_{X,K}\right\|_{\Ec}^2
\sum_{k=0}^{N-1} \frac{1}{2M_k} \leq   \sum_{k=0}^{N-1} \frac{1}{2M_k}\left\|\tg_{X,k}\right\|_{\Ec}^2  \leq \psi(x_0) - \psi(x_{N}) + \frac{\e}{4} \sum_{k=0}^{N-1} \frac{1}{M_k} + N(4\delta_u + \delta_{pu}).
\notag
\end{align}
Finally, since, for all $x\in X$ $\psi(x) \geq \psi^* > -\infty$ and $\tg_{X,K}\stackrel{\eqref{eq:gXeqMxmx}} {=} M_K (x_K - x_{K+1})$, we obtain
\begin{align}
&\left\|M_K (x_K - x_{K+1}))\right\|_{\Ec}^2  \leq \left(\sum_{k=0}^{N-1} \frac{1}{2M_k} \right)^{-1} (\psi(x_0) - \psi^* + N(4\delta_u + \delta_{pu})) + \frac{\e}{2},
\end{align}
which is \eqref{eq:pg_la_2_rate}.
The estimate for the number of checks of Inequality~\eqref{eq:UGMCheck} is proved in the same way as in \cite{nesterov2006cubic}, but we provide the proof for the reader's convenience.
Let $i_k \geq 1$ be the total number of checks of  Inequality~\eqref{eq:UGMCheck} on the step $k \geq 0$. Then $i_0 = 1+\log_2\frac{M_0}{L_0}$ and, for $k \geq 1$, $M_k = 2^{i_k-1}L_k = 2^{i_k-1}\frac{M_{k-1}}{2}$. Thus, $i_k = 2+ \log_2\frac{M_k}{M_{k-1}}$, $k \geq 1$. Then, the total number of checks of Inequality~\eqref{eq:UGMCheck} is 
\begin{equation}
\sum_{k=0}^{N-1}i_k=1+\log_2\frac{M_0}{L_0} + \sum_{k=1}^{N-1}\left(2+ \log_2\frac{M_k}{M_{k-1}}\right) = 2N-1+\log_2\frac{M_{N-1}}{L_0}.
\end{equation}
\qed

Let us consider two corollaries of the theorem above. First is a simple case, when in Definition \ref{D:IO} $L(\delta_c) \equiv L$. Second is the case, when $L(\delta_c)$ is given by \eqref{eq:Lofd}.
\begin{Cor}
\label{Cor:Sm}
Assume that there exists a constant $L>0$ s.t. for the dependence $L(\delta_c)$ in Definition \ref{D:IO} it holds that $L(\delta_c) \leq L$ for all $\delta_c > 0$. Assume also that there exists a number $\psi^* > -\infty$ such that $\psi(x) \geq \psi^*$ for all $x \in X$. Then, after $N$ iterations of Algorithm \ref{Alg:UGM}, it holds that
\begin{equation}
\left\|M_K (x_K - x_{K+1}))\right\|_{\Ec}^2  \leq \frac{4L(\psi(x_0) - \psi^*)}{N} + 4L(4\delta_u + \delta_{pu}) + \frac{\e}{2}.
\label{eq:pg_la_2_rate_1}
\end{equation}
Moreover, the total number of checks of Inequality \eqref{eq:UGMCheck} is not more than 
$$
2N+\log_2\frac{L}{L_0}.
$$
\end{Cor}
\myproof
By our assumptions, for all iterations $k\geq 0$, there exists $i_k \geq 0$ s.t. $M_k = 2^{i_k}L_k \geq L(\delta_{c,k}) \equiv L$. Hence, we can apply Theorem \ref{Th:UGMRate}. Let $i_k \geq 1$ be the total number of checks of Inequality \eqref{eq:UGMCheck} on a step $k\geq 0$. Then, for all $k\geq 0$, the inequality $M_k = 2^{i_k}L_k \leq 2L$ should hold. Otherwise the termination of the inner cycle would happen earlier. Using this inequalities, we obtain
$$
\left(\sum_{k=0}^{N-1} \frac{1}{2M_k} \right)^{-1} \leq \left(\sum_{k=0}^{N-1} \frac{1}{4L} \right)^{-1} = \frac{4L}{N}.
$$
Thus \eqref{eq:pg_la_2_rate_1} follows from Theorem \ref{Th:UGMRate}. The same argument proves the second statement of the corollary.
\qed
\begin{Cor}
\label{Cor:Hold}
Assume that the dependence $L(\delta_c)$ in Definition \ref{D:IO} is given by \eqref{eq:Lofd} for some $\nu \in (0,1]$, i.e. 
\begin{equation}
L(\delta_c) = \left( \frac{1-\nu}{1+\nu} \cdot \frac{2}{\delta_c} \right)^{\frac{1-\nu}{1+\nu}} L_{\nu}^{\frac{2}{1+\nu}}, \quad \delta_c > 0.
\label{eq:Lofd1}
\end{equation}
Assume also that there exists a number $\psi^* > -\infty$ such that $\psi(x) \geq \psi^*$ for all $x \in X$. 
Then, after $N$ iterations of Algorithm \ref{Alg:UGM}, it holds that
\begin{equation}
\left\|M_K (x_K - x_{K+1}))\right\|_{\Ec}^2  \leq 2^{\frac{1+3\nu}{2\nu}} \left( \frac{1-\nu}{1+\nu} \cdot \frac{40}{\e} \right)^{\frac{1-\nu}{2\nu}} L_{\nu}^{\frac{1}{\nu}} \left(\frac{\psi(x_0) - \psi^*}{N} + 
(4\delta_u + \delta_{pu}) \right) + \frac{\e}{2}.
\label{eq:pg_la_2_rate_2}
\end{equation}
Moreover, the total number of checks of Inequality \eqref{eq:UGMCheck} is not more than 
$$
2N-1+\frac{1+\nu}{2\nu}+\frac{1-\nu}{2\nu}\log_2\left(40 \cdot \frac{1-\nu}{1+\nu}  \right) + \frac{1-\nu}{2\nu} \log_2\frac{1}{\e} +  \log_2 \frac{L_{\nu}^{\frac{1}{\nu}}}{L_0}.
$$
\end{Cor}
\myproof
First, let us check that, for any constants $c_1, c_2 >0$, there exists an integer $i \geq 0$ s.t. $2^ic_1 \geq L\left(\frac{c_2}{c_12^i}\right)$. Substituting $\delta_c = \frac{c_2}{c_12^i}$ to \eqref{eq:Lofd1} gives 
$$
L\left(\frac{c_2}{c_12^i}\right) = 2^{\frac{1-\nu}{1+\nu}i} c_3,
$$
where $c_3>0$ is some constant. Since $1-\frac{1-\nu}{1+\nu} = \frac{2\nu}{1+\nu} >0$, we conclude that the required $i \geq 0$ exists. Thus, we can apply Theorem \ref{Th:UGMRate}.

Let $i_k \geq 1$ be the total number of checks of Inequality \eqref{eq:UGMCheck} on a step $k\geq 0$. Then, for all $k\geq 0$, the inequality $M_k = 2^{i_k}L_k \leq 2L(\delta_{c,k})$ should hold. Otherwise the termination of the inner cycle would happen earlier. 
From this inequality and \eqref{eq:Lofd1} it follows that
\begin{equation}
M_k \leq 2 \left( \frac{1-\nu}{1+\nu} \cdot \frac{40M_k}{\e} \right)^{\frac{1-\nu}{1+\nu}} L_{\nu}^{\frac{2}{1+\nu}}.
\end{equation}
Solving this inequality for $M_k$, we obtain
\begin{equation}
M_k \leq  2^{\frac{1+\nu}{2\nu}} \left( \frac{1-\nu}{1+\nu} \cdot \frac{40}{\e} \right)^{\frac{1-\nu}{2\nu}} L_{\nu}^{\frac{1}{\nu}}.
\label{eq:MkBound}
\end{equation}
Whence,
\begin{equation}
\left(\sum_{k=0}^{N-1} \frac{1}{2M_k} \right)^{-1} \leq \left(\sum_{k=0}^{N-1} \frac{1}{4L} \right)^{-1} = 2^{\frac{1+3\nu}{2\nu}} \left( \frac{1-\nu}{1+\nu} \cdot \frac{40}{\e} \right)^{\frac{1-\nu}{2\nu}} \frac{L_{\nu}^{\frac{1}{\nu}}}{N}.
\end{equation}
Now \eqref{eq:pg_la_2_rate_2} follows from Theorem \ref{Th:UGMRate}.

Using \eqref{eq:InnerChecksBound} and the bound \eqref{eq:MkBound}, we obtain the estimate for the total number of checks of Inequality \eqref{eq:UGMCheck}.
\qed

Let us make some remarks about the obtained results. First, if we set in Corollary \ref{Cor:Hold} $\nu=1$, we recover the result of Corollary \ref{Cor:Sm}. 
Second, in the situation of Corollary \ref{Cor:Hold}, to make the controlled part of the right-hand side smaller than $\e$ we need to choose 
$$
N \geq {\rm const} \cdot \frac{L_{\nu}^{\frac{1}{\nu}}(\psi(x_0)-\psi^*)}{\e^{\frac{1+\nu}{2\nu}}}.
$$
One can see that the less $\nu$ is, the worse is the bound. This is expected as for non-smooth non-convex problems the norm of gradient mapping $g_X(\cdot)$ at the stationary point could not be equal to zero.
Third, we can see that uncontrolled error $4\delta_u + \delta_{pu}$ can dramatically influence the error estimate, especially, when $\nu$ tends to zero. 

Finally, let us explain, why small $\left\|M_K (x_K - x_{K+1}))\right\|_{\Ec}$ means that $x_{K+1}$ is a good approximation for stationary point of the initial problem \eqref{eq:PrStateInit}. Let us prove the following result, which was communicated to us by Prof. Yu. Nesterov without proof.
\begin{Lm}
\label{Lm:NecOptCond}
Let in Problem \eqref{eq:PrStateInit} $f(x)$ be continuously differentiable, $h(x)$ be convex, $X$ be a closed convex set. Assume that $x^*$ is a local minimum in this problem.
Then, for all $x \in X$,
\begin{equation}
\la \nabla f(x^*), x- x^* \ra + h(x) + h(x^*) \geq 0.
\label{eq:NecOptCond}
\end{equation} 
\end{Lm}
\myproof Let us fix an arbitrary point $x \in X$. Denote $x_t = tx+(1-t)x^* \in X$, $t \in [0,1]$. Since $x^*$ is a local minimum in \eqref{eq:PrStateInit}, $X$ is a convex set, $h(x)$ is a convex function, we obtain for all sufficiently small $t > 0$
$$
0 \leq \frac{f(x_t)+h(x_t)-f(x^*)-h(x^*)}{t} \leq \frac{f(x_t)-f(x^*)}{t} + h(x)-h(x^*).
$$
Taking the limit as $t \to +0$, we prove the stated inequality.

Assume, for simplicity, that we are in the situation of Subsection \ref{SS:SFIO}. This means that $f(x)$ is $L(f)$-smooth, we can uniformly approximate its gradient
\begin{equation}
\|\bar{g}(x)-\nabla f(x)\|_{{\Ec},*} \leq \bar{\delta}_c^2+\bar{\delta}_u^2,
\label{eq:nfappr}
\end{equation}
and the set $X$ is bounded with diameter $D$. Also assume that the chosen prox-function $d(\cdot)$ is $L(d)$-smooth.

From \eqref{eq:InPrMap}, \eqref{eq:InPrMap1}, \eqref{eq:UGMwStep}, we obtain that there exists $\nabla h(x_{K+1}) \in \partial h(x_{K+1})$ s.t., for all $x \in X$,
$$
\left\la \tg(x_K,\delta_{c,K},\delta_u) + M_K\left[d'(x_{K+1}) - d'(x_K) \right]  + \nabla h(x_{K+1}), x - x_{K+1}\right\ra \geq -\delta_{pc,K} - \delta_{pu}. 
$$
Whence, by convexity of $h(x)$,
\begin{align}
\la \nabla f(x_{K+1}) , x- x_{K+1} \ra + h(x) - h(x_{K+1}) \geq &\la \nabla f(x_{K+1}) -\nabla f(x_{K}) , x- x_{K+1} \ra \\
&+ \la \nabla f(x_{K}) - \tg(x_k,\delta_{c,k},\delta_u) , x- x_{K+1} \ra \\
&+ \la M_k\left[d'(x_{K}) - d'(x_{K+1}) \right]  , x - x_{K+1}\ra -\delta_{pc,K} - \delta_{pu}, \quad x \in X.
\label{eq:LongOptCondIneq}
\end{align}
By $L(f)$-smoothness of $f$, boundedness of $X$, we obtain
$$
\la \nabla f(x_{K+1}) -\nabla f(x_{K}) , x- x_{K+1} \ra \geq - \frac{L(f)}{M_K} \|M_K(x_{K} - x_{K+1})\|_{\Ec}D.
$$
From \eqref{eq:nfappr}, by boundedness of $X$, we get
$$
 \la \nabla f(x_{K}) - \tg(x_K,\delta_{c,K},\delta_u) , x- x_{K+1} \ra \geq - (\bar{\delta}_{c,K}^2+\bar{\delta}_u^2) D.
$$
Using $L(d)$ smoothness of $d(x)$ and boundedness of $X$, we obtain
$$
\la M_k\left[d'(x_{K}) - d'(x_{K+1}) \right]  , x - x_{K+1}\ra \geq -L(d) \|M_K(x_{K} - x_{K+1})\|_{\Ec}D.
$$
Substituting last three inequalities to \eqref{eq:LongOptCondIneq}, we obtain that, if $\|M_K(x_{K} - x_{K+1})\|_{\Ec} \leq \e$, then
$$
\la \nabla f(x_{K+1}) , x- x_{K+1} \ra + h(x) - h(x_{K+1}) \geq -\Theta(\e) - \bar{\delta}_u^2 D - \delta_{pu}.
$$
Thus, at the point $x_{K+1}$ the necessary condition in Lemma \ref{Lm:NecOptCond} approximately holds.
%
%
%
%
%

\section*{Conclusion}
In this article, we propose a new adaptive gradient method for non-convex composite optimization problems with inexact oracle and inexact proximal mapping. We showed that, for problems with inexact H\"older-continuous gradient, our method is universal in terms of H\"older parameter and constant. For the proposed method, we prove convergence theorem in terms of generalized gradient mapping and show that a point returned by our algorithm is a point where necessary optimality condition approximately holds.

\textbf{Acknowledgments.} The author is very grateful to Prof. A. Nemirovski, Prof. Yu. Nesterov, Prof. B. Polyak for fruitful discussions.

\bibliographystyle{plainnat}
\bibliography{references}


\end{document}